\documentclass[graybox]{svmult}

\usepackage[mathcal]{euscript}
\usepackage[all]{xy}
\usepackage{color,soul} 
\usepackage{hyperref} 
\hypersetup{
 colorlinks=true, 	 
 linktoc=all,     	 
 linkcolor=black,  	 
 citecolor=blue,    
 filecolor=black,  
 urlcolor=black       
}

\usepackage{mathptmx}       
\usepackage{helvet}         
\usepackage{courier}        
\usepackage{type1cm}        
\usepackage{makeidx}         
\usepackage{graphicx}        
\usepackage{multicol}        
\usepackage[bottom]{footmisc}

\usepackage{amscd,amsmath,amssymb,mathrsfs}

\makeindex             

\newcommand{\R}{{\mathbb R}}
\newcommand{\C}{{\mathbb C}}

\newcommand{\T}{{\mathcal T}}

\begin{document}

\title*{\Large A Survey on Spaces of Homomorphisms to Lie Groups}
\titlerunning{A Survey on Spaces of Homomorphisms to Lie Groups} 

\author{Frederick R. Cohen and Mentor Stafa}
\authorrunning{F. R. Cohen and M. Stafa} 

\institute{
Frederick R. Cohen 
\at University of Rochester, Rochester NY 14627, USA
\at\email{fred.cohen@rochester.edu}
\and 
Mentor Stafa 
\at Tulane University, New Orleans LA 70118, USA
\at {\it Current address}: ETH Z\"urich, Switzerland
\at\email{mstafa@tulane.edu}
}

%
%

\maketitle

\vspace{-1in}

\begin{flushright}
\it Dedicated to Samuel Gitler Hammer,\\ who brought us much joy and interest in Mathematics.
\end{flushright}

\abstract*{ The purpose of this article is to give an exposition of  topological properties of 
spaces of homomorphisms from certain finitely generated discrete groups to Lie groups $G$, 
and to describe their connections to classical representation theory,  as well as other structures. 
Various properties are given when $G$ is replaced by a small category, or the discrete 
group is given by a right-angled Artin group.}

\abstract{The purpose of this article is to give an exposition of  topological properties 
of spaces of homomorphisms from certain finitely generated discrete groups to Lie 
groups $G$, and to describe their connections to classical representation theory,  
as well as other structures. Various properties are given when $G$ is replaced by a 
small category, or the discrete group is given by a right-angled Artin group.}

\keywords{space of homomorphisms, representation variety, small category, right-angled Artin group}


\section{Introduction}

Let $\pi$ be a finitely generated discrete group and $G$ a Lie group. In this paper we 
describe features of spaces of group homomorphisms from $\pi$ to $G$ denoted by 
$Hom(\pi,G)$, which is given the subspace topology of the Cartesian product of a 
finite number of copies of $G$.

Spaces of group homomorphisms $Hom(\pi,G)$
exhibit interesting properties. They have played an important role in mathematical work 
related to physics going back to work of E. Witten \cite{witten1,witten2}, who considered 
spaces of commuting pairs and triples in certain Lie groups. Following the same lines but 
using different methods A. Borel, R. Friedman and J. Morgan developed the space of 
commuting  as well as \textit{almost commuting} pairs and triples \cite{borel2002almost}, 
i.e. pairs and triples commuting up to an element in the center. The spaces have also 
been studied by V. Kac and A. Smilga \cite{kac.smilga} and W. Goldman \cite{goldman} 
as well as T. Friedmann and R. Stanley \cite{friedmann2014counting}.

The group $G = G^{ad}$ acts by conjugation on the space $Hom(\pi,G)$. The quotient 
space $\mathcal{M}(\pi,g)= Hom(\pi,G)/G^{ad}$ is called the \textit{representation space} 
(or \textit{moduli space}  or \textit{representation variety}, \textit{character variety}) 
of $\pi$ in $G$. Representation varieties 
for fundamental groups of surfaces in compact Lie groups have been studied in the 
context of gauge theories and hyperbolic geometry when $G$ is the complex or real 
projective special linear group. These spaces also admit natural symplectic  structures 
as discovered by W. Goldman \cite{goldman1984symplectic}.

Recently the subject experienced mathematical development in a different direction. 
Here, we restrict to homotopy theoretic properties of $Hom(\pi,G)$. A. Adem and 
F. Cohen \cite{adem2007commuting} considered the problem in the framework 
of spaces of homomorphisms $Hom(\pi,G)$, where $\pi$ is the free abelian group 
of any rank $n$. These spaces split when suspended once as given in Theorem 
\ref{thm: splitting of Hom(Zn,G)}.

Properties of $Hom(\pi,G)$ have been developed in \cite{friedmann2014counting} 
in case $G$ is a finite group.  Using a different approach, the authors 
\cite{stafa.thesis,stafa.comm} assemble all the spaces $Hom(\mathbb{Z}^n,G)$ 
into a larger, more tractable space $Comm(G)$ reminiscent of a Stiefel manifold. 
This is discussed in more detail in the following sections. A few problems at the 
interface of these subjects are posed at the end of this article.

The spaces $Hom(\mathbb{Z}^n,G)$ also admit the structure of simplicial spaces. 
Their geometric realization, denoted by $B(2,G)$ was introduced by A. Adem, F. Cohen 
and E. Torres-Giese \cite{fredb2g}.
One curious connection arises by turning the natural map 
$$B(2,G) \to BG$$ 
into a fibration with homotopy fibre denoted $E(2,G)$.  Basic properties of the 
monodromy representation of $G$ on the first homology group of the fibre are curious. 
For example, if $G$ is a finite group of odd order, then the classical Feit-Thompson 
theorem is equivalent to the map $H_1(E(2,G))_ G  \to H_1(B(2,G)),$ where $H_1(E(2,G))_ G$ 
denotes the module of coinvariants under the monodromy action, failing to be surjective. 
Properties of the monodromy representation may be interesting, but are currently 
not well understood \cite{fredb2g,stafa.monodromy}.
This setting suggests that the spaces $B(2,G)$ contain compelling information. 
Thus a categorical interpretation is given here in section \ref{sec: B(q,G)}.

Much beautiful work has been done by many people, some of whom are not mentioned here because of the brevity of this article. Some of this work is surveyed below.


\section{Spaces of homomorphisms and their topology}

Let $\pi$ be a finitely generated discrete group with $n$ generators $x_1,\dots,x_n$. 
Then a homomorphism $f\in Hom(\pi,G)$ is determined by the images of the generators 
$x_1,\dots,x_n$ in $G$. Thus there is the natural identification of a homomorphism  
$f$ with the $n$-tuple $(f(x_i),\dots,f(x_n))$ using the natural inclusion
$$
Hom(\pi,G) \hookrightarrow G^n.
$$
The space of homomorphisms can be endowed with the subspace topology of $G^n$. 
Equivalently, this association of the subspace topology can be seen from a choice of 
surjection from the free group with $n$ letters to $\pi$
$$
F_n \to  \pi,
$$
which induces an injection on the level of spaces
$$
Hom(\pi,G) \hookrightarrow Hom(F_n,G)=G^n.
$$
In particular, the image of this inclusion is homeomorphic to $
Hom(\pi,G) $ in case $G$ is a Lie group.
The space $Hom(\pi,G)$ is regarded as a pointed space with base-point given by 
the constant function $f_0$ with $f_0(x)=1_G \in G$, for all $x\in \pi$.


\section{Path-components and fundamental groups}

V.~Kac and A.~Smilga \cite{kac.smilga} addressed the number of path-components 
in the case of $Hom(\mathbb{Z}^n,G)$ for $G = Spin(7)$. The number of path 
components for the case when $G=SO(k), \ \mbox{or}  \  O(k)$ was calculated by 
G. H. Rojo \cite{rojo}. He shows that if $k$ is large enough, the number of path 
components of $Hom(\mathbb{Z}^n,G)$ no longer depends on $k$.

Applying the \textit{classifying space functor} gives a map
$$ 
B: Hom(\pi,G) \rightarrow Map_{\ast}(B\pi,BG),
$$
for which each group homomorphism $f$ gives a pointed map $Bf:B\pi \to BG$.
It should be noted that if $G$ is compactly generated as a topological space, then 
the classifying space functor is a continuous functor, thus $B$ is continuous 
\cite{Steenrod}.

Applying $\pi_0$, the path components, there is an induced map $B_0$ 
$$B_0: \pi_0(Hom(\pi,G)) \rightarrow [B\pi,BG],$$
where the homotopy classes of maps $[B\pi,BG]$ is in one-to-one correspondence 
with all principal $G$-bundles over $B\pi$. Hence, properties of $\pi_0(Hom(\pi,G))$ 
inform on principal $G$-bundles over $B\pi$, even though the map $B_0$ is neither 
an injection nor a surjection, in general.

Moreover, the space of homomorphisms can exhibit curious properties such as 
being disconnected even if   $G$ is 1-connected. The case for which every abelian 
subgroup of a compact and connected Lie group $G$ can be conjugated to a 
subgroup of the maximal torus $T$ of $G$ is special.
In this case, the space of ordered commuting $n$-tuples $Hom(\mathbb{Z}^n,G)$ 
is path-connected as given in \cite[Proposition 2.3]{adem2007commuting}.

A related natural conjecture was that if $G$ is a finite, discrete group, then $B(2,G)$ 
is a $K(\pi,1)$ based on examples in \cite{fredb2g}. C.~Okay in his thesis  \cite{okay} 
gave a counter-example in the special case where $G$ is an extra-special $2$-group of order $32$.

More recent work of  S. Lawton and D. Ramras \cite{ramras} addressed the fundamental 
group of the character variety of $Hom(\pi,G)/G^{ad}$ where $\pi$ is the fundamental 
group of a compact orientable surface. Tom Baird studied the individual spaces 
$Hom(\mathbb Z^n,G)$ via the Leray spectral sequence of the inclusion 
$Hom(\mathbb Z^n,G) \hookrightarrow G^n$, where $G$ is a compact Lie group 
\cite{bairdcohomology}.  Baird gave the additive structure for the rational cohomology 
of  $Hom(\mathbb Z^n,G)$, which is also an isomorphism of graded rings, in case 
$G$ is compact, and connected.


\section{Commuting elements for non-compact $G$}

Suppose that $\pi=\mathbb{Z}^n$ and that $G$ is a reductive algebraic group. 
A. Pettet and J. Souto \cite{pettet.souto} showed that if $K$ is a maximal compact 
subgroup of a reductive algebraic group $G$, then the inclusion 
$$Hom(\mathbb{Z}^n,K)\subset Hom(\mathbb{Z}^n,G)$$ 
is a strong deformation retract. This result has been generalized by M. Bergeron 
\cite{bergeron} to finitely generated nilpotent groups. One of the corollaries of 
\cite{pettet.souto} is the isomorphism 
$$ \pi_1(Hom(\mathbb{Z}^n,G))\cong \pi_1(G)^n,$$ 
which was first proven in \cite{gomez.pettet.souto}. Therefore we will usually 
restrict attention to $G$ given by a compact and connected Lie group.

The above feature is special for the case of  $Hom(\mathbb{Z}^n,G)$. Namely, 
if $K$ is a maximal compact subgroup of $G$, the natural map  $Hom(\pi,K) \to  Hom(\pi,G)$, 
where $K$ is the maximal compact subgroup of $G$, may not induce an isomorphism 
on the set of path-components.
Many such examples arise with $G = PSL(2,\R)$, and $K = SO(2)$ where $\pi$ is 
given by the fundamental group of a closed, orientable surface of genus $g > 1$. 
Other examples arise by setting $G = PSL(2,\C)$ with $K = SO(3)$ and $\pi$ equal 
to certain mapping class groups \cite{adem2007commuting}.


\section{Stable decompositions}\label{section:stable}

A general method of stabe decompositions was given in [1]. It was shown in  
\cite[Theorem 1.6]{adem2007commuting} that the space $Hom(\mathbb{Z}^n,G)$ 
admits a stable decomposition as a wedge sum of certain spaces resembling a 
fat wedge with  appropriate definitions given next.

\begin{definition}\label{defn: singular subspace)}
Define $S_k(G)\subset Hom(\mathbb{Z}^k,G)$ to be the subspace of ordered 
$n$-tuples  having at least one coordinate equal to the identity element.  
Define $$\widehat{Hom}(\mathbb{Z}^k,G) = Hom(\mathbb{Z}^k,G)/S_k(G).$$
\end{definition}

It was observed in \cite{adem2007commuting} that the inclusion 
$S_k(G)\subset Hom(\mathbb{Z}^k,G)$  is a cofibration in case $G$ is  a Lie group. 
Thus there is a cofibre sequence 
$$S_k(G) \to Hom(\mathbb{Z}^k,G) \to  \widehat{Hom}(\mathbb{Z}^k,G).$$

This basic point-set topological property is useful in what follows.
\begin{theorem}[{{\cite[Theorem 1.6]{adem2007commuting}}}]\label{thm: splitting of Hom(Zn,G)}
Assume that $G$ is a closed subgroup of $GL_m(\C)$. Then (i) the cofibre sequence 
$S_k(G) \to Hom(\mathbb{Z}^k,G) \to  \widehat{Hom}(\mathbb{Z}^k,G)$ 
after suspending once is split, and (ii)
there is a homotopy equivalence
$$
\Sigma\left( Hom(\mathbb{Z}^n,G) \right) \to \bigvee_{1 \leq k \leq n} \Sigma \left( \bigvee_{n \choose k} \widehat{Hom}(\mathbb{Z}^k,G) \right).
$$

\end{theorem}

This decomposition follows from the natural decomposition of $Hom(A \times B,G)$ 
where $A$, and $B$ are discrete groups. This observation is addressed in 
Lemma \ref{lem: free products} below.

The spaces $Hom(\mathbb{Z}^n,G)$ can be assembled to form a new space 
called $Comm(G)$ as a subspace inside the free associative monoid generated by $G$, 
the \textit{James reduced product} $J(G)$ of $G$ \cite{james1955reduced}, as follows.

\begin{definition}[{{\cite[Definition 1.1]{stafa.thesis}}}]\label{defn: Comm(G)}
Let $G$ be a topological group. Define
$$
Comm(G)=\left( \bigsqcup_{k \geq 0} Hom(\mathbb{Z}^k,G) \right)/\sim
$$
where (i) $Hom(\mathbb Z^0, G) = \{1_G\}$, and (ii) the equivalence relation $\sim$ 
is generated by the relation
$$
(g_1,\dots,g_k)\sim (g_1,\dots,\widehat{g_i},\dots,g_k) \text{ if } g_i=1_G.
$$
\end{definition}

One feature of the space $Comm(G)$ is that it can be regarded as a universal construction 
which contains all of the spaces $Hom(\mathbb{Z}^k,G) $, for all $k$, and is also 
computationally tractable. Namely, homological properties of the individual spaces 
$Hom(\mathbb{Z}^n,G)$ are obtained directly from $Comm(G)$. 
That follows as $Comm(G)$ admits a stable decomposition which directly implies 
information about $Hom(\mathbb{Z}^k,G)$ for all $k$. In addition, an equivariant 
function space is given below which gives the 
homology groups of $Comm(G)$.

\begin{theorem}[{{\cite[Theorem 1.12]{stafa.thesis}}}]\label{thm: stable decomp of Comm(G)}
Let $G$ be a Lie group. Then there is a homotopy equivalence
$$
\Sigma(Comm(G)) \to \bigvee_{n \geq 1} \Sigma(\widehat{Hom}(\mathbb{Z}^n,G)).
$$
\end{theorem}

The statement of Theorem \ref{thm: stable decomp of Comm(G)} is that the suspension 
of $Comm(G)$ contains exactly one copy of the {suspension of the } quotient
$$\widehat{Hom}(\mathbb{Z}^n,G)=Hom(\mathbb{Z}^n,G)/S_n(G)$$ 
that appears in the decomposition of $\Sigma(Hom(\mathbb{Z}^n,G))$, for each integer $n > 0$.

\begin{corollary}\label{thm: stable decomp of any homology of Comm(G)}
Let $G$ be a Lie group, and let $E_*(-)$ be any reduced homology theory.
Then there are isomorphisms
\begin{enumerate}
	\item $E_*(Comm(G)) \to \oplus_{1 \leq j < \infty} E_*(\widehat{Hom}(\mathbb{Z}^j,G))$, and
	\item $E_*(Hom(\mathbb Z^n,  G)) \to \oplus_{1 \leq j \leq n} \oplus_{\binom{n}{j}}E_*(\widehat{Hom}(\mathbb{Z}^j,G))$.
\end{enumerate}
\end{corollary}

Next assume that $G$ is a compact and connected Lie group with maximal torus $T$, 
and Weyl group $W$ such that $Comm(G)$ is path-connected. Then the homology 
groups of $Comm(G)$ are organized into a relatively simple construction.
Let $\mathcal T[V]$ denote the tensor algebra generated by the module $V$ where 
all modules are assumed to be free over the ring $R = \mathbb Z[1/|W|]$.
With these conditions on $G$, for which the order of the Weyl group is a unit in 
{\it ungraded} singular reduced homology $H_*(-;R)$, then the following holds 
where $UH_*(X;R)$ denotes {\it ungraded homology}.

\begin{corollary}\label{thm: stable ungraded homology for Comm(G)}
Let $G$ be a compact and connected Lie group with maximal torus $T$, 
and Weyl group $W$ such that $Comm(G)$ is path-connected.
Consider {\it ungraded} singular reduced homology $UH_*(-;R)$ where $R$ is the ring
$R = \mathbb Z[1/|W|]$. Then the ungraded homology $UH_*(Comm(G);R)$ is 
a free module over $R = \mathbb Z[1/|W|]$ which is isomorphic to 
$\mathcal T(\widetilde{H}_*(T; R))$.
\end{corollary}

This corollary follows from Theorem \ref{thm: homology of Comm(G)} and is 
restated more generally in Section \ref{sec: homological results} 
(see also \cite[Theorem 1.17]{stafa.comm}). Some examples are the exceptional 
Lie groups $G_2,\, F_4,\, E_6,\, E_7,\, E_8$ and $SU(n),\,U(n),\,Sp(n),\,Spin(n)$, 
where the ungraded rational homology of $Comm(G)$ is isomorphic to a tensor 
algebra generated by a module of rank equal to the rank of the reduced 
rational homology of the maximal torus.

A more precise hold on $Comm(G)$ will be discussed in the next few sections. 
For example, the next theorem exhibits an equivalence from an equivariant 
function space to $Comm(G)$.

\begin{theorem}\label{thm: stable homology for Comm(G)}
If $G$ is compact, simply-connected Lie group, then there is a map 
$$G\times_{NT} \Omega\Sigma (T) \to Comm(G),$$ 
which induces an equivalence in homology after the order of the Weyl group $W$ has been inverted. 
Thus the cohomology of $Comm(G)$ with the order of $W$ inverted is the ring of $W$-invariants
in the cohomology of $$G/T \times \Omega\Sigma (T).$$
\end{theorem}
 
T.~ Baird gave the first computation of the rational homology of 
$Hom(\mathbb{Z}^k,G)$ in case $G$ is simply-connected, and compact 
in \cite{bairdcohomology}. D.~Sjerve, and E.~Torres-Giese
gave (i) a complete computation for the homology of
$Hom(\mathbb{Z}^k,G)$ in case $G = SU(2)$, and (ii) they also gave a 
stable decomposition in terms of Thom spaces of multiples of the classical 
Hopf line bundle over $\mathbb R \mathbb P^2$ \cite{Sjerve-Torres}. 
T.~Baird, L.Jeffrey, and P.~Selick gave an independent analogous result \cite{bjs}.


\section{Related constructions associated to $G$}

Suppose that $\pi= F_n/\Gamma^q$ where $\Gamma^q$ denotes the $q$-th 
stage of the descending central series for $F_n$. Then there is a generalization 
of $Comm(G)$ as follows.

\begin{definition}[{{\cite[Definition 1.1]{stafa.thesis}}}]\label{defn: X(,q,G)}
Let $G$ be a topological group. Define
$$
X(q,G)=\left( \bigsqcup_{k \geq 1} Hom(F_n/\Gamma^q,G) \right)/\sim
$$
where the relation $\sim$ is generated by the relation
$$
(g_1,\dots,g_k)\sim (g_1,\dots,\widehat{g_i},\dots,g_k) \text{ if } g_i=1_G.
$$
\end{definition}

Note that $X(2,G)=Comm(G)$. A natural filtration of $J(G)$ is obtained 
$$
X(2,G)\subset X(3,G)\subset \cdots \subset X(q,G)\subset \cdots \subset J(G).
$$
There are analogues where (i) the descending central series is replaced by 
the mod-$p$ descending central series, or (ii) $F_n$ is replaced by either 
the pro-finite completion, or pro-$p$ completion. Similar results are satisfied, 
but not addressed here.

Little is known about the features of $X(q,G)$ at the moment other than the 
existence of a similar stable decomposition as in the case of $Comm(G)$.
Namely, each such space of homomorphisms $Hom(F_n/\Gamma^q,G) $ again 
splits after suspending into (1) a ``singular summand", and (2) a complementary part. 
The next Theorem was stated, but without full details of proof in 
\cite{adem2007commuting}. It was proven  by B. Villarreal (private communication) 
noting that the "singular set" $S_k(G) \subset Hom(F_n/\Gamma^q,G) $ is a closed 
subspace of a real algebraic set which is triangulable.

\begin{theorem}\label{thm: splitting of Hom(Fn/q,G)}
Assume $G$ is a closed subgroup of $GL_m(\C)$. Then there is a homotopy equivalence
$$
\Sigma\left( Hom(F_n/\Gamma^q,G) \right) \to \bigvee_{1 \leq k \leq n} 
		\Sigma \left( \bigvee_{n \choose k} Hom(F_k/\Gamma^q,G)/S_k(G)\right).
$$ 
 where $S_k(G)\subset Hom(F_k/\Gamma^q,G)$ is equal to the subset with 
 elements having at least one coordinate the identity element.
\end{theorem}

The next decomposition was given in \cite{stafa.thesis,stafa.comm}:

\begin{theorem}[{{\cite{stafa.thesis,stafa.comm}}}]\label{thm: stable decomp of X(q,G)}
Let $G$ be a Lie group. Then there is a homotopy equivalence
$$
\Sigma(X(q,G)) \to \bigvee_{n \geq 1} \Sigma(Hom(F_n/\Gamma^q,G)/S_n(G)).
$$
\end{theorem}

One compelling reason for considering these spaces is that their homology provides interesting
representations of $Aut(F_n/\Gamma^q)$ especially interesting in case $G$ is a finite group.
Related problems are considered in the problem section here.


\section{Approximations of $Comm(G)$}\label{sec: approximations of Comm(G)}

For any abelian subgroup $ A \subset G$, there is a map
\begin{align*}
G\times A^n &\to Hom(\mathbb{Z}^n,G)\\
(g,g_1,\dots ,g_n) & \mapsto (g_1^g,\dots ,g_n^g)
\end{align*}
gotten by conjugating the $n$-tuple $(g_1,\dots ,g_n)$ by $g$. This map factors 
through the quotient obtained from the natural action of the  normalizer 
$NA$ of $A$ to give a map
$$
G\times_{NA}A^n \to Hom(\mathbb{Z}^n,G).
$$
Finally, applying the construction of Definition \ref{defn: Comm(G)} one obtains a map
$$
\alpha: G\times_{NA}J(A) \to Comm(G),
$$
which is central to our study of the spaces $Hom(\mathbb{Z}^n,G)$ for the 
case when $A$ is the maximal torus $T$ of $G$. Let $W$ be the \textit{Weyl group} 
of $G$ defined as the quotient $W=NT/T$. Then the map $\alpha$ can be rewritten as 
$$
\alpha: G/T\times_{W}J(T) \to Comm(G).
$$

Now let $Hom(\mathbb{Z}^n,G)_1$ denote the path component of the trivial 
homomorphism. Define $Comm(G)_1$ to be the space 
$$
Comm(G)_1=\left(\bigsqcup_{n\geq 1} Hom(\mathbb{Z}^n,G)_1\right)/\sim
$$
with the same relation. Then the map $\alpha$ restricted to $Comm(G)_1$ is a surjection. 
The space $G\times_{NT}J(T)$ also admits  stable decomposition of an equivariant 
function space \cite[Theorem 1.9]{stafa.thesis} as follows.

\begin{theorem}[{{\cite[Theorem 1.13]{stafa.thesis}}}]\label{thm: stable decomp Gx_NT J(T)}
Let $G$ be a Lie group. Then there are homotopy equivalences
$$
\Sigma\big(G\times_{NT}J(T)\big) \to \Sigma\bigg(G/T \vee 
	\bigvee_{q\geq 1}(G\times_{NT}\widehat{T}^q)/(G/NT) \bigg),
$$
where $\widehat{T}^q$ is the $q$-fold smash product of $T$.
\end{theorem}

With theorems \ref{thm: stable decomp of Comm(G)} and \ref{thm: stable decomp Gx_NT J(T)} 
one can show the following.

\begin{theorem}[{{ \cite[Theorem 1.13]{stafa.thesis}}}]\label{thm: alpha Gx_NT J(T) to Comm(G)}
Let $G$ be a compact Lie group. Then the map 
$$\alpha: G/T\times_{W}J(T) \to Comm(G)_1$$ 
induces an isomorphism in homology if the order of the Weyl group is inverted.  
\end{theorem}

Thus $Comm(G)$ can be approximated by simultaneously conjugating products 
of the maximal torus $T$, equivalently by conjugating elements of $J(T)$. 
This description of $Comm(G)$ and the approximation allow for a direct 
investigation of the spaces $Hom(\mathbb{Z}^n,G)$  individually. In particular, 
the additive homology of these spaces can be obtained this way as 
shown in the next section.


\section{Homological results}\label{sec: homological results}

Recall that the map $\alpha: G/T\times_{W}J(T) \to Comm(G)_1$ induces an 
equivalence in homology if the order of the Weyl group is inverted. Therefore 
one obtains the homology of $Comm(G)_1$ with coefficients in 
$R=\mathbb{Z}[|W|^{-1}]$ from the homology of $G/T\times_{W}J(T)$, 
where $\mathbb{Z}[|W|^{-1}]$ is the ring of integers with the order of the 
Weyl group inverted. Let $R[W]$ denote the group ring, and let $V$ be the 
reduced homology of $T$ as an $R[W]$-module with $\T[V]$ the tensor algebra 
of $V$. Using the Leray-Serre spectral sequence for homology, the homology 
of $Comm(G)_1$ with coefficients in $R$ is obtained as follows.

\begin{theorem}\label{thm: homology of Comm(G)}
Let $G$ be a compact and connected Lie group with maximal torus $T$ and 
Weyl group $W$. Then there is an isomorphism in homology
$$
H_{\ast}(Comm(G)_1;R)\cong H_{\ast}(G/T;R)\otimes_{R[W]}\T[V].
$$
\end{theorem}

Note that in Theorem \ref{thm: homology of Comm(G)}, if the grading of 
homology is not considered, namely, ungraded homology is taken, where 
the order of the Weyl group is inverted, then the ungraded homology of the 
flag variety $G/T$, denoted $UH_*(G/T)$, is isomorphic as an $R[W]$-module 
to the group ring $RW$ itself \cite{bairdcohomology}. 
Let $UH_{\ast}(X;R)$ denote the \textit{ungraded homology} of the space $X$, 
the usual singular homology, but regraded with $UH_*(X;R) = \oplus_{0\leq j} H_j(X;R)$. 
Since $Comm(G)_1$ is not necessarily equal to $Comm(G)$, 
Corollary \ref{thm: stable ungraded homology for Comm(G)} is a special case 
of the following immediate corollary of Theorem \ref{thm: homology of Comm(G)}.

\begin{corollary}
There is an isomorphism in ungraded homology given by
$$
UH_{\ast}(Comm(G)_1;R)\cong \T[V].
$$
\end{corollary}

Thus the homology groups of $Comm(G)$ with real coefficients are, roughly speaking, 
reassembling the tensor algebra by regrading the tensor algebra.

\begin{definition}\label{defn: trigrading)}
Consider cohomology of $Comm(G)$ with real coefficients. A tri-graded Poincar\'e series 
$$\sum_{i,j,m}A(i,j,m)q^i s^j t^m$$ 
is defined for the real cohomology of $Comm(G)_1$ obtained from the module of invariants 
$$\bigg(H^{\ast}(G/T;\R)\otimes \T^{\ast}[V]\bigg)^W.$$ Here $\T^{\ast}$ 
denotes the dual of the tensor algebra, and the coefficient $A(i,j,m)$ is the rank 
of the invariants of the module obtained by tensoring the $i$-degree cohomology 
of $G/T$ with the dual of an $m$-fold tensor in $H^{\ast}(J(T);\R) = \T
[V]$ in homological degree $j$. 
Equivalently, this module $H^{\ast}(G/T;\R)\otimes \T^{\ast}[V]$ is given as a direct sum of modules
$$
M_{(i,j,m)}=\sum_{\substack{j=k_1+\cdots+k_m \\ k_q>0}} 
	H^i \otimes (\Lambda^{k_1}\R^n \otimes \cdots \otimes \Lambda^{k_m}\R^n).
$$
\end{definition}

Consider cohomology with real coefficients.  Denote the module 
$H_{\ast}(G/T;\R)$ by $H$. An application of Molien's theorem gives 
the coefficients $A(i,j,m)$ in the Poincar\'e series with the following closed 
form as worked out by Vic Reiner in the appendix to \cite{stafa.comm}.

\begin{theorem}[{{Reiner \cite[Theorem 1.20]{stafa.comm}}}]\label{thm: Poincare series for Comm(G)}
If $G$ is a compact and connected Lie group with maximal torus $T$ and 
Weyl group $W$, then the tri-graded Poincar\'e series of $Comm(G)$ is
$$
\begin{aligned}
&P((H\otimes \T^{\ast}[V])^W;q,s,t) \\
&\hspace*{.5in} =\frac{\prod_{i=1}^n (1-q^{2d_i})}{|W|} 
	\sum_{w \in W} \frac{1}{\det(1-q^2w)(1-t(\det(1+sw)-1))}.
\end{aligned}
$$
\end{theorem}

A corollary of this theorem gives the reduced homology of all spaces 
$Hom(\mathbb{Z}^n,G)$ with coefficients in $R$.

\begin{corollary}
Let $G$ be a compact and connected Lie group. Then the reduced cohomology 
of $Hom(\mathbb{Z}^m,G)$ is given additively by
$$
\widetilde{H}^d(Hom(\mathbb{Z}^m,G);R)\cong \sum_{1\leq s \leq m} 
		\sum_{i+j=d}\left( \sum_{\substack{j=k_1+\cdots+k_s \\ i\geq 0}} 
		\bigoplus_{m \choose s}\left( M_{(i,j,s)}\right)^W\right).
$$
\end{corollary}


\section{The spaces $B(q,G)$}\label{sec: B(q,G)}

A special case of natural subspaces of the classifying space $BG$ were defined in \cite{fredb2g} using
$Hom(\mathbb{Z}^m,G)$ as follows.

J.~Milgram's construction of $BG$ is the geometric realization of a simplicial space 
which has objects in degree $n$ given by $G^n$ \cite{milgram1967bar}.
The subspaces $Hom(\mathbb{Z}^n,G) \subset G^n$ are preserved by the face and 
degeneracy operations. Thus there is an associated geometric realization denoted 
$B(2,G)$ \cite{fredb2g}.

Similar constructions apply in the case of $Hom(F_n/\Gamma^q, G)$ to obtain a space $B(q,G)$.
There are inclusions $$B(2,G) \subset B(3,G) \subset \cdots \subset BG.$$ Consider 
the homotopy theoretic fibre $B(q,G) \to BG$ to obtain a fibration $$E(q,G) \to B(q,G) \to BG.$$ 
The Serre exact sequence for the fibration $$E(q,G) \to B(q,G) \to BG$$ has several natural properties. 
Namely, consider the associated  Serre exact sequence 
$$H_1(E(q,G))_G \to H_1(B(q,G)) \to H_1(BG) \to 0.$$  
If $G$ is finite and odd order, then the map $H_1(E(q,G))_G \to H_1(B(q,G))$  is not a surjection. 
However the action of $G$ on $H_1(E(q,G))$ is still not well-understood, and an 
independent proof of the previous fact is not known.  One problem is 
listed below concerning this action.

In addition, recall that Milnor showed that if $X$ is a path-connected CW-complex, 
then there is a topological group $G(X)$ such that $BG(X)$ has the homotopy type 
of $X$ \cite{milnor56b}. Thus the spaces $B(q,G(X))$ give a new filtration of any 
connected CW-complex $X$
$$B(2,G(X)) \to B(3,G(X)) \to \cdots \to BG(X) \simeq X.$$ 
Properties of this filtration as well as other related filtrations were introduced in 
\cite{fredb2g}. One immediate corollary follows next, see also 
\cite[Theorem 6.3]{fredb2g} for a proof.

\begin{corollary}\label{cor: B(q,G)}
Let $G(X)$ denote a topological group whose classifying space is homotopy 
equivalent to $X$, a path-connected finite CW-complex. Then the looped 
maps $\Omega(B(q,G(X))) \to \Omega(BG(X))$
all have a section up to homotopy.
\end{corollary}

The spaces $B(2,G)$ and their generalizations have several further properties described here.
The cohomology of $B(2,G)$ with the order of the Weyl group inverted 
was described in \cite{fredb2g}. Examples were given by transitively 
commutative groups, and Suzuki groups.

Adem and Gomez \cite[Definition 2.1]{adem.gomez3} give the next definition.

\begin{definition}\label{defn: stuff)}
Assume that $X$ is a CW-complex. A principal $G$-bundle $q:E \to X$ is 
{\it transitionally commutative} provided there is an open cover $\{U_i \ |  \ i \in I \}$ of $X$ 
such that the bundle $q: E \to X$ restricted to $U_i$ is trivial, and the 
transition functions $$\rho(i,j) : U_i \cap U_j \to G$$ commute when they 
are defined simultaneously.
\end{definition}

\begin{theorem}[{\cite[Theorem 2.2]{adem.gomez3}}]
Suppose that G is a Lie group. Let $f : X \to BG$ denote the classifying map of a 
principal G-bundle $q : E \to X$ over the finite CW-complex X . Then up to homotopy, 
$f$ factors through $B(2,G)$ if and only if there is an open cover of X on which the 
bundle is trivial over each open set and such that on intersections the transition 
functions commute when they are simultaneously defined, i.e. $q$ is transitionally commutative.
\end{theorem}

A beautiful theorem of Adem, Gomez, Lind, Tillmann \cite{adem.gomez.lind.tillmann} follows next.

\begin{theorem}[{\cite[Theorem 1.1]{adem.gomez.lind.tillmann}}]
The spaces $B(q,G )$ for $$G = O, SO, U, SU, Sp$$ provide a filtration by 
$E_{\infty}$-ring spaces of the classical infinite loop spaces BSU, BU, BSO, BO and 
BSp, respectively. 
\end{theorem}

Recall the natural maps $$\Sigma_n \to O(n)$$ induce an infinite loop map 
$$\Omega^{\infty} S^{\infty} \to \mathbb Z \times  BO$$ closely tied to the 
J-homomorphism. Since this map may be chosen to be a homomorphism, there is 
an induced map $$B(q,\Omega^{\infty}S^{\infty}) \to B(q,\mathbb Z \times BO).$$ 
It is reasonable to ask about the behavior of this map. 

The purpose of the rest of  this section is to give an extension of the construction 
$B(2,G)$ where $G$ is replaced by a small category $\mathcal C$ to obtain a space 
denoted $$B(2,\mathcal C).$$ The methods are based on the article by 
M.~Weiss \cite{Weiss}. This extension is also a direct consequence of \cite{adem.gomez3}.

The classifying space of a topological group $G$ is the geometric realization of the 
simplicial space gotten from all ordered $n$-tuples of composable morphisms. 
The space $B(2,G)$ is analogous where all elements in an ordered $n$-tuple are 
required to pairwise commute. Thus given any topological category $\mathcal C$, 
there is an analogous classifying space obtained from ordered $n$-tuples of 
composable elements which themselves pairwise commute. Furthermore, 
these morphisms are endomorphisms of the same object. The associated 
classifying space is denoted $B(2, \mathcal C).$

 Weiss considers the contravariant functors from $\mathcal C$ to the category of 
 sets together with natural transformations between these functors. The space 
 $B\mathcal C$ classifies representable sheaves over a topological space with 
 stalks which are commutative diagram of composable endomorphisms. In the 
 process of proving this, Weiss refines a $\mathcal C$-set as contravariant functors 
 from $\mathcal C$ to sets with the morphisms given by natural transformations. 
 The $\mathcal C$-sets which are of the form $b \to mor_{\mathcal C}(b,c)$ 
 are said to be representable.

There is a sheaf theoretic analogue of the definition stated next which  Adem 
and Gomez gave in \cite[Definition 2.1]{adem.gomez3}.

\begin{definition}\label{defn: transitionally commutative sheaves)}
Assume that $X$ is a CW-complex, and $\mathcal C$ is a small category.
A $\mathcal C$-sheaf over $q:\mathcal {E} \to X$ is {\it transitionally commutative} 
provided there is an open cover of $X$ $\{U_i \ |  \ i \in I \}$
such that the sheaf $q: \mathcal{E} \to X$ restricted to $U_i$ is trivial, and the transition 
functions $$\rho(i,j) : U_i \cap U_j \to q^{-1}(s \in U_i\cap U_j)$$ commute when they 
are defined simultaneously.
\end{definition}

\begin{theorem}\label{thm: B(2,Cat)}
If $\mathcal C$ is a small category, the classifying space $B(2, \mathcal C)$ classifies representable
sheaves of $\mathcal C$-sets over a topological space which are transitionally commutative.
\end{theorem}

The proof of this is essentially that given in \cite{Weiss, adem.gomez3}.


\section{Right-angled Artin groups}\label{section:RAAGs}

The purpose of this section is to demonstrate how earlier methods imply 
stable decompositions for $Hom(\pi,G)$, where $\pi$ is a finitely 
generated right-angled Artin group.

Recall that a finitely generated right-angled Artin group is given by the fundamental group
of the polyhedral product $Z(K;(S^1, *))$, where $K$ is a finite simplicial complex \cite{bbcg,davis.okun,stafa.monodromy}.  In particular, a presentation is given 
by generators $x_1, \dots, x_m$ with relations $[x_i,x_j] = 1$ if and only if vertices 
$i$ and $j$ share an edge in $K$.

Write $\pi(K)$ for the fundamental group of  $Z(K;(S^1, *))$.
The fundamental groups of $Z(K;(S^1, *))$ are all right-angled Artin groups. 
The group $\pi(K)$ is given by a graph product of groups \cite{davis.okun, stafa.monodromy}. In
addition, if $K$ is a flag complex, then  $Z(K;(S^1, *))$ is a $K(\pi,1)$  \cite{davis.okun, stafa.monodromy}.

\begin{definition}
Let $\sigma$ denote a {\it maximal full face} with $m(\sigma)$ vertices. Namely, 
(i) $\sigma$ is a face which is abstractly a $(j-1)$-simplex $\Delta[j-1]$ with $j= m(\sigma)$ 
vertices as a simplicial complex, and (ii) if $\sigma$ is contained in any face $\tau$ in $K$, 
then $\sigma = \tau$.  Let $\pi(\sigma)$ denote the fundamental group of
$Z(\sigma;(S^1, *)) = (S^1)^j$. Thus $\pi(\sigma)$ is a free abelian group of rank 
$m(\sigma)$  the number of vertices of $\sigma$.
Let $\mathcal{M}(K)$ denote the set of maximal full faces of $K$.
\end{definition}

The next Lemma stated without proof is an observation, where $\ast_{\alpha} G_{\alpha}$ 
denotes the free product of groups $G_{\alpha}$.

\begin{lemma}\label{lem: }
Let $\pi(K)$ denote the fundamental group of the polyhedral product $Z(K;(S^1, *))$.
Then there are surjections
$$\ast_{\sigma \in \mathcal{M}(K) }  \pi(\sigma) \to \pi(K),$$ and
$$\prod_{\sigma  \in \mathcal{M}(K)   }  \pi(\sigma) \to \mathbb Z^m.$$
\end{lemma}

The next Lemma gives general properties about $Hom(\pi,G)$.

\begin{lemma}\label{lem: free products}
Let $G$ be a Lie group, and $A,B,C$ be discrete groups. Then
\begin{enumerate}
\item The natural forgetful map $$Hom(A \ast B,G) \to Hom(A,G) \times Hom(B,G)$$ is a homeomorphism.
\item If $f:A \to C$ is a surjection, then the induced map $$Hom(C,G) \to Hom (A,G)$$ is an embedding.
\item The group $Aut(A) \times Aut(G)$ acts naturally on $Hom(A,G)$.
\item If $G$ is simply-connected and compact, then $Hom(\mathbb Z^n, G)$ is a closed subset of $G^n$.
\item The suspension $\Sigma(Hom(A \times B,G))$ is homotopy equivalent to 
$$\Sigma(  (Hom(A,G) \vee   (Hom(B,G)) \vee W(A,B, G)) $$
where $W(A,B,G)$ depends on $A,B, G.$
\end{enumerate}
\end{lemma}

One immediate consequence follows.

\begin{corollary}\label{cor: free products}
Let $G$ be a Lie group, and $A_i = \mathbb Z^{n_i}$, $ 1 \leq i \leq m$. Then the space
$Hom(  {\ast}_{1\leq i \leq m}A_{i}, G) $ is homeomorphic to $   
\prod_{1 \leq i \leq m} Hom( {\mathbb Z}^{n_i},G).$
\end{corollary}

\begin{corollary}\label{cor: splitting}
Let $G$ be a Lie group.
Let $K$ be a finite simplicial complex, with  $\pi(K)$ the fundamental group of 
$Z(K;(S^1, *))$, and $\mathcal{M}(K)$ the set of maximal faces of $K$. Then $$Hom(\pi(K), G) \subset     
\prod_{\sigma \in \mathcal{M}(K)} Hom(\pi(\sigma), G)$$ where $\pi(\sigma)$ is isomorphic to
$\mathbb \mathbb{Z}^{m(\sigma)}$.
Furthermore, there is a homotopy equivalence
$$\Sigma(Hom(\pi(K), G)) \to \Sigma(X)  \vee \bigvee_{\sigma \in \mathcal{M}(K)} \Sigma(Hom(\mathbb{\mathbb{Z}}^{m(\sigma)}, G))$$ for a space $X =  Hom(\pi(\sigma), G)/ \bigvee_{\sigma \in \mathcal{M}(K)}(Hom(\mathbb{\mathbb{Z}}^{m(\sigma)}, G)) $. 
\end{corollary}

\begin{example}\label{exm: splitting}
Let $K$ be a flag complex with $6$ vertices $\{1,2,3,4,5,6\}$ and with maximal faces
$\sigma =\{1,2,3\}$, $\tau = \{4,5,6\}$, and $ \rho = \{3,4\}$. Then $Hom(\pi(K),G)$ 
has the following three spaces as a retract:  $Hom(\pi(\sigma),G)$, 
$Hom(\pi(\tau),G)$, and $Hom(\pi(\rho),G)$
It follows by a direct computation that there is a homotopy equivalence
$$\Sigma(Hom(\pi(K),G)) \to \Sigma X \vee \Sigma(Hom(\pi(\sigma),G)) \vee 
	\Sigma(Hom(\pi(\tau),G)) \vee \Sigma(Hom(\pi(\rho),G))$$ 
for some space $X$. Furthermore, 
$$Hom(\pi(K),G) =  \bigg{ (  }   Hom(\pi(\sigma),G)) \times Hom(\pi(\tau),G)) \bigg{)}      
	\cap \bigg{(} G^2 \times   Hom(\pi(\rho),G)  \times G^2  \bigg{)}.$$ 
\end{example}

\begin{definition}
Let $\mathcal{M}(K)$ denote the set of maximal full faces of $K$ which has $m$ vertices.
Let $\sigma$ denote a maximal full face in $\mathcal{M}(K)$. Define 
$$G(\sigma)  = \{(g_1,g_2, \cdots, g_m) \in G^m |  \  [g_i,g_j] = 1 \mbox{ for all}  \ i,j \in \sigma\}.$$
\end{definition}

Observe that $$Hom(\pi(K), G) \subset    G(\sigma)$$  for each $ \sigma \in \mathcal{M}(K)$. 
Thus $Hom(\pi(K), G) \subset \bigg{ (  }     \bigcap_{\sigma \in \mathcal{M}(K)   }   G(\sigma) \bigg{ )  }.$ 
Also, observe that $  \bigcap_{\sigma \in \mathcal{M}(K)   }   G(\sigma) \subset Hom(\pi(K), G)$ 
restated as a Theorem.

\begin{theorem}\label{thm: RAAG decomposition}
If $G$ is a compact Lie group, then
$$Hom(\pi(K), G) =  \bigcap_{\sigma \in \mathcal{M}(K)   }   G(\sigma) .$$ 
where, $G(\sigma)$ is homeomorphic to $Hom(\mathbb Z^t,G) \times G^{m-t}$ for some $t$
by a homeomorphism which permutes coordinates, and $Hom(\mathbb Z^t,G)$
is a retract of $Hom(\pi(K), G)$ for each $\sigma \in \mathcal{M}(K)$,  where $\sigma$ has $t$ vertices.
\end{theorem}

\begin{corollary}\label{thm: MVSS for RAAGS}
If $G$ is a compact Lie group, then $$G^m - Hom(\pi(K), G) = \bigcup_{\sigma \in \mathcal{M}(K)   } \big(   G^m - G(\sigma)\big)$$ gives a Mayer-Vietoris spectral sequence abutting to the homology of \newline
$G^m - Hom(\pi(K), G)$.
\end{corollary}


\section{Problems}\label{sec: problems}
This section is a list of certain problems concerning the spaces constructed here.

\begin{enumerate}

\item  Let $G$ be a Lie group and $\mathcal{C}(G)$ be the space of closed 
subgroups of $G$ as defined in \cite{Chabauty}, and developed in 
Bridson-de la Harpe-Kleptsyn \cite{Bridson.Harpe.Kleptsyn}.

\begin{definition}\label{defn: Chabauty topology)}
The \textit{Chabauty topology} on $\mathcal{C}(G)$ is defined by giving a 
basis of neighborhoods for a closed subgroup $H \in \mathcal C(G)$ by the subsets
$$
\mathcal{V}_{K,U}(H)=\{D\in \mathcal{C}(G): D \cap K \subset HU \text{ and } H \cap K \subset DU\},
$$
where $K$ is compact in $G$ and $U$ is an open neighborhood of the identity $1 \in G$. 
A subspace of $\mathcal{C}(G)$ which is of interest is the 
\textit{space of abelian closed subgroups} $\mathcal{A}(G)$. 
\end{definition}

Describe the analogue of the Chabauty topology for closed, and maximal abelian subgroups. 
Does $Comm(G)$ admit a map onto this space? If there is such a map, is it a quasi-fibration?

\item Recall that an infinite loop space has the homotopy type of a topological group.
In view of the results of  \cite{adem.gomez.lind.tillmann}, it is natural to ask whether 
the spaces $B(q,G)$ are infinite loop spaces in more generality.

Do infinite loop maps induce infinite loops on the level of $B(q,G)$?

Replace $G$ by a small category $\mathcal C$. Does $B(q,G)$ admit a natural 
extension to $B(q, \mathcal C)$? If so, what are the properties?

Is the homotopy type of $B(q,G)$ independent
of the multiplication in $G$ where $G$ is an infinite loop space ? Note that
$BG$ is not an invariant of the homotopy type of $G$
in general.

\item Any reasonable space $X$  is the classifying space of a group $G(X)$. 
What does $B(q, G(X))$ look like?  This construction thus gives a new filtration of $X$. 
It is natural to ask whether this filtration impacts the space $X$ in a useful way. 
For example, is this filtration a homotopy invariant?



\item Fix a discrete group $G$, and ask which representations of $Aut(\pi)$ 
occur in the natural action on $Hom(\pi,G)$. If $G$ is finite, and $\pi$ is finitely 
generated, this process gives permutation representations of $Aut(\pi)$.

An example is given next from unpublished work of Shiu-Chun Wong, Jie Wu 
and the first author given by families of  two-stage nilpotent finite groups 
$\Gamma$ which are quotients of $Aut(F_n)$, but where the quotient map 
does not factor through $GL(n,\mathbb Z)$ for $n >1$.

For instance, the quotient of the free group of rank 2 by the third stage of 
its descending central series, denoted $F_2/\Gamma^3$, 
is a quotient of $Aut(F_2)$ by the usual action. 
Requiring the generators of  $F_2/\Gamma^3$, say $x, y$, to have order 
$p$ as well as $[x,y]$  to have order $p$, gives a group $G$ of order $p^3$.
Checking Magnus' well-known presentation 
for $Aut(F_2)$, there is an induced epimorphism
$Aut(F_2) \to G$
where $G$ is the group of order $p^3$ defined earlier,
for $p> 3$.
Thus $G$ is a quotient of $Aut(F_2)$ by a check of the classical presentation, 
but not a quotient of $GL(2,\mathbb{Z})$.

 Which finite groups are quotients of $Aut(F_2)$ obtained from the natural action 
 on the finite  set $Hom(F_2,G)$ for $G$ finite? It is classical that 
 $PGL(2,\mathbb F_q)$ is such a quotient where $\mathbb F_q$ denotes a finite field.

Does anything ``new" occur when considering  the action of 
$Aut(F_2/\Gamma^q)$ on $Hom(F_2/\Gamma^q,G)$?
 
 Let $\hat{F_n}$ denote the pro-finite completion of the free group 
 on $n$ letters. Let $G$ denote a finite group. Thus $Aut(\hat{F_n})$ 
 acts naturally on $Hom(\hat{F_n}, G).$ Are all finite groups given by quotients of this action?
Which finite groups occur as quotients of $Aut(\hat{F}_2)$ acting on the set
$Hom(\hat{F_2}, G)$ in this setting? Observe that there are finite groups 
which are quotients  of $Aut(F_2)$, but not quotients of $GL(2,\mathbb Z)$.

There are natural representations of $Aut(\hat{F_n})$  on the homology of the space
$Hom(\hat{F_n}, G).$  What are these representations?




\item Give information about the monodromy representations for $$B(q,G) \to BG.$$ Show that
the map $H_1(E(q,G))_G \to H_1B(q,G)$ is not onto if $G$ is of odd order.
\end{enumerate}

\noindent{\bf Acknowledgements}  {\small The authors would like to thank the 
\textit{Centro di Ricerca Matematica Ennio De Giorgi} at the Scuola Normale Superiore 
in Pisa, Italy, and the \textit{Istituto Nazionale di Alta Matematica} for their support. 
Frederick R. Cohen was partially supported by the Institute of Mathematics and 
its Applications. 
Mentor Stafa was partially supported by DARPA grant number N66001-11-1-4132.}



\end{document}